\documentclass[12pt]{article}
\usepackage{amsmath,amsfonts,amssymb,amsthm}
 \usepackage{authblk}

\title{On the rate of Poisson approximation \\to partial sum processes\footnote{This work was supported by the  Russian Foundation for Basic Research, under grant
\newline 18--01--00074.}}

\date{}

\author[1,2]{Pavel S. Ruzankin\thanks{ruzankin@math.nsc.ru}}
\author[1,2]{Igor S. Borisov\thanks{sibam@math.nsc.ru}}
\affil[1]{\normalsize Sobolev Institute of Mathematics, Novosibirsk, Russia}
\affil[2]{Novosibirsk State University, Novosibirsk, Russia}

\begin{document}

\maketitle

 \begin{abstract}
We investigate approximation of a Bernoulli partial sum process to the accompanying Poisson process in the non-i.i.d. case.
The rate of closeness is studied in terms of the minimal distance in probability.
\end{abstract}

\section{Introduction and results}

Let $X,X_1,X_2,\dots$~
be independent random variables taking values in a measurable Banach space
$\big(\bf  B,\cal A,\|\cdot \|\big)$,
where $\cal A$ is a $\sigma$-field of subsets of ${\bf  B}$.
 Denote by $P, P_1,P_2,\dots,$ the respective distributions.
We will assume that the linear operations and the norm
$\|\cdot\|$
are ${\cal A}$-measurable.
As a rule, if ${\bf  B}$ is separable, then ${\cal A}$ is
the Borel $\sigma$-field.
The above-mentioned measurability conditions are fulfilled in this case.
For nonseparable Banach spaces, the class $\cal A$ is usually smaller
than the Borel $\sigma$-field (for example, $\cal A$ may be the cylindrical or the ball
$\sigma$-field).

Denote by $Pois(\mu)$ the compound Poisson distribution with  L\'{e}vy
measure $\mu$:
$$
Pois(\mu):=e^{-\mu({\bf  B})}\sum_{k=0}^{\infty}\frac{\mu^{*k}}{k!},
$$
where $\mu^{*k}$~ is the $k$-fold convolution of a finite measure $\mu$
with itself; $\mu^{*0}$ is the unit mass at zero.
As is easy to see, $Pois(\mu)$ coincides with
the distribution of the sum
$\sum\nolimits_{i\le\pi(\mu({\bf  B}))}Y_i$,
where $\{Y_i\}$ are independent identically  distributed  random variables
with the distribution  $\mu(\cdot)/\mu({\bf  B})$, and
 $\pi (\mu({\bf  B}) )$~
is a  Poisson random variable with mean $\mu({\bf  B})$, which is independent
of $\{Y_i\}$. (Here $\sum\nolimits_{i\le 0}=0$ by definition.)

The main subject of the paper is the two partial sum processes
\begin{equation}
S_n\equiv S_n(t):= \sum_{i\le nt}X_i, \,\,\,\,t\in [0,1],
                                                        \label{eqno1}
\end{equation}
and
\begin{equation}
\Pi_n\equiv \Pi_n(t):= \sum_{i\le nt}\hat{\pi}_i(P_i), \,\,\,\,t\in [0,1],
                                                             \label{eqno2}
\end{equation}
where $ \{\hat{\pi}_k(P_k) \}$~
are independent random variables with the respective distributions
$ \{Pois(P_k) \}$. We consider these ${\bf B}$-valued stochastic processes as elements
of the Banach space of ${\bf B}$-valued right-continuous functions on $[0,1]$ endowed with the sup-norm.
We study the proximity in terms of the distance
$$d(z,S_n,\Pi_n):=\inf_{\{X_k,\hat{\pi}_k(\cdot);\ k\le n\}}{\bf  P}
 \biggl(\sup_{t\in [0,1]}\|S_n(t)-\Pi_n(t)\|>z \biggr)
$$
\begin{equation}
=\inf_{\{X_k,\hat{\pi}_k(\cdot);\ k\le n\}}{\bf  P}
 \biggl(\max_{k\le n}\biggl\|\sum_{i=1}^{k}(X_i-\hat{\pi}_i(P_i))
\biggr\|>z \biggr)
\label{eqno3}
\end{equation}
where the infimum is taken over all collections of
$\{X_k,\,i\le n\}$ and $\{\hat{\pi}_k(\cdot),\,i\le n\}$
defined on a common probability space. It is worth noting that this distance belongs to the class of so-called minimal distances
between distribu\-tions of ${\bf B}$-valued stochastic processes with paths in the above-mentioned  Banach space
(e.g., see Zolotarev, 1976).

Alongside with (\ref{eqno3}), we also need the
total variation
distance between the distribu\-tions
${\cal L}(Y_i)$
of arbitrary ${\bf  B}$-valued random variables $Y_i$, $i=1,2$:
$$
V (Y_1, Y_2  ):=
\sup_{A\in {\cal A}} |{\cal L}(Y_1)(A)-{\cal L}(Y_2)(A)|.
$$
There is the following connection between these two distances:
$$V(S_n,\Pi_n)=d(0,S_n,\Pi_n).
$$
This equality is a consequence of the following duality theorem by Dob\-rushin (1970):
\begin{equation}
V (X, Y)=\inf_{X,Y}{\bf  P}(X\ne Y),
\label{eqno4}
\end{equation}
where the infimum
is taken over all pairs
($X, Y$) based on a common probability space.
Relation (\ref{eqno4}) means that if the total variation distance between some random
variables is sufficiently small, then there exist versions of these random
variables coinciding with probability close to $1$.

Systematic investigation of the accuracy of Poisson approximation to the distri\-butions of the random variable $S_n(\ref{eqno1})$ as well as of the stochastic process $S_n(\cdot)$ began in the middle of 20-th century. We would like to note the pioneer paper by Prohorov (1953) and the three papers by Le Cam (1960, 1965, 1970)
where unimprovable estimates (up an absolute constant factor) for the total variation distance $V(S_n(\ref{eqno1}),\Pi_n(\ref{eqno1}))$ were obtained.
Since that time, a significant number of published studies were devoted to this problem.
Regarding estimation of the total variation distance we refer the reader to the references in the book by  Barbour et al. (1992) and the recent review by Novak (2019b).
Notice that there is a number of papers devoted to the Poisson approximation in terms of other probability distances like $d(z,\cdot,\cdot)$
(e.g., see Ruzankin, 2001), the information divergence (e.g., see Harremo\"{e}s and Ruzankin, 2004), $\chi^2$-distance (e.g., see Borisov and Vorozheikin, 2008),
and of other distances (e.g., see Novak, 2019a; Ruzankin, 2004, 2010).

We would like to  recall the most significant  relevant results.
In the case
${\bf  B}={\bf  R}$, let $X_i\equiv\nu_i$~
be Bernoulli random variables with the respective success probabilities $p_i$. Denote by $\pi_i$,\, $i=1,\ldots,n$, independent
Poissonian random variables with respective means $p_i, \,i=1,\ldots,n$.
There is the following two-sided estimate obtained by Barbour and Hall (1984):
\begin{equation}
\frac{1}{32}
\varepsilon
\sum_{i\le n}p_i^2\le V
\Bigl(\sum_{i\le n}\nu_i, \sum_{i\le n}\pi_i\Bigr)\Bigr)
\le
\varepsilon
\sum_{i\le n}p_i^2,
\label{eqno5}
\end{equation}
where
$\varepsilon:=\min\left\{ 1,\,\left(\sum_{i\le n}p_i\right)^{-1}\right\}$.
It is worth noting that, in general, for non-Bernoulli sequences $\{X_i\}$ of nonidentically distributed random variables,
upper and lower bounds for the total variation distance
will be slightly worse than (\ref{eqno5}) and has the unimprovable order $O\left(\sum_{i\le n}p_i^2\right)$ (see Le Cam, 1960; Borisov, 1996). 

To formulate an analog of inequality (\ref{eqno5}) for non-Bernoulli sequences $\{X_i\}$ we need the notation
$$\overline{\nu}_n(k):=\sum\nolimits_{i\le k}\nu_i,\,\,\,\,
\bar{\pi}_n(k):=\sum\nolimits_{i\le k}\pi_i,\,\,\,
k\le n.
$$
Denote by
$\{X_i^0\}$
 independent random variables with the distributions
\begin{equation}
{\cal L}(X_i^0)={\cal L}(X_i\mid X_i\ne 0).
                                                   \label{eqno6}
\end{equation}
Note that these random variables are well defined because, by the above assump\-tions,
the event $\{X_i\ne 0\}= \{\|X_i\|\ne 0 \}$  is measurable.
It is proved in Borisov (1996)
that, for every $k$, the following equalities hold:
\begin{equation}
{\cal L} (X_k) ={\cal L}
\biggl(\,\sum\limits_{i=\overline{\nu}_n(k-1)+1}^{\overline{\nu}_n(k)}
X_{k,i}^0\biggr),\quad
Pois (P_k)={\cal L}
\biggl(\,\sum\limits_{i=\bar{\pi}_n(k-1)+1}^{\bar{\pi}_n(k)}
X_{k,i}^0\biggr),
                                                            \label{eqno7}
\end{equation}
where $\bigl\{X_{k,i}^0; i=1,2,\dots\bigr\}$~are independent copies of
$X_k^0$,
$k\le n$, which do not depend on the random processes $\overline{\nu}_n(\cdot)$ and $\overline{\pi}_n(\cdot)$, with $p_k={\bf P}(\|X_k\|\ne 0)$, $k\le n$.
Notice that the first relation in (\ref{eqno7}) can be rewritten as ${\cal L} (X_k)={\cal L} (\nu_kX^0_k)$ but the second one cannot.
If $X_i^0$ are {\it identically distributed}
(whereas $X_i$ may be nonidentically
distributed with arbitrary $p_i$),
then (\ref{eqno13}) and independence of the
increments of the processes $\bar{\nu}_n(\cdot)$ and $\bar{\pi}_n(\cdot)$
readily imply more informative Khintchine's representations  for the
distributions of
$S_n$ and $\Pi_n$ in the Banach space $\bf  B^n$:
\begin{equation}
{\cal L}(S_n)={\cal L}\left\{\sum\limits_{i=1}^{\bar{\nu}_n(k)}
X_i^0;\ k=1,\dots,n\right\},\quad
{\cal L}(\Pi_n)={\cal L}\left\{\sum\limits_{i=1}^{\bar{\pi}_n(k)}
X_i^0;\ k=1,\dots,n\right\}.                                                          \label{eqno8}
\end{equation}
The relations (\ref{eqno4}), (\ref{eqno7}) and (\ref{eqno8}) imply the following results (see Borisov, 1996):

{\bf Lemma 1}. {\it For arbitrary distributions $\{P_i\}$,
$$
V (S_n, \Pi_n)\le V (\bar{\nu}_n(\cdot), \bar{\pi}_n(\cdot) ).
                                                                  \label{eqno9}
$$
}

{\bf Lemma 2}. {\it
If
${\bf  P} \bigl(\sum\nolimits_{i\le m}X_{k,i}^0=0\bigr)=0$
for all natural numbers $m$ and $k$ then
\begin{equation}
V (S_n, \Pi_n  )=V (\bar{\nu}_n(\cdot),
\bar{\pi}_n(\cdot) ).
                                             \label{eqno10}
\end{equation}
}

{\bf Lemma 3}. {\it
If
$\|X_i\|\le 1$  with probability $1$ for all $i\ge 1$, and the random variables $\{X_i^0\}$
are identically distributed then, for all $z\ge 0$,
\begin{equation}
d(z,S_n,\Pi_n)\le d (z,\bar{\nu}_n(\cdot),\bar{\pi}_n(\cdot)).
                                                                        \label{eqno11}
\end{equation}
}

Finally, the
last assertion is connected with a special structure of
 $\{X_i\}$. We will say that $X_i$ {\it are elements of indicator type
with conforming supports} if, for all natural $m$ and $k$, the following identity
is valid with probability 1:
\begin{equation}
 \left\|\sum\limits_{i=m+1}^{m+k}X_i^0\right\|=k.
\label{eqno12}
\end{equation}

{\bf Lemma 4}. {\it
If $\{X_i\}$ are random elements of indicator type with conforming supports
and the random variables $\bigl\{X_i^0\bigr\}$ are identically distributed then, for all $z\ge 0$,
\begin{equation}
d(z,S_n,\Pi_n)=d (z,\bar{\nu}_n(\cdot),\bar{\pi}_n(\cdot)).
\label{eqno13}
\end{equation}
}

It is worth noting that, to study the generic problem of Poisson approximation, one can find one or another implicit form of reducing the problem
to that for the Bernoulli case  (e.g., see Borisov and Ruzankin, 2002).
The Lemmas above allow us to do it the shortest way.

In the case of  i.i.d. random variables $\{X_i\}$ one has

{\bf Theorem 1} (Borisov, 1996). {\it
Let $\{\nu_i\}$~ be independent identically distributed Bernoulli random variables
with success probability $p$. Then, for all $z\ge 0$,
$$
d (z,\bar{\nu}_n(\cdot),\bar{\pi}_n(\cdot))\le
(np^2)^{\lfloor z\rfloor+1}\exp \{-C_1z\ln\ln(z+2)+C_2\}\,\,\,\,
\text{if}\,\,\,  np\ge 1,$$
\begin{equation}
d (z,\bar{\nu}_n(\cdot),\bar{\pi}_n(\cdot))\le np^{\lfloor z\rfloor+2}\exp\{-C_3z+C_4\}\,\,\,\,\text{if} \,\,\, np\le 1,
                                                               \label{eqno14}
\end{equation}
where $C_1$--\,$C_4$~ are absolute positive constants.
}

Adduce an example of random elements of indicator type with conforming supports.
 Consider the problem of uniform Poisson approximation
of so-called local multivariate empirical processes. Let
$F_n(\bf  z)$, ${\bf z}\in{\bf  R}^k,$ be the empirical distribution
function based on the sample $y_1,\dots,y_n$ of i.i.d. vectors
with an arbitrary distribution in ${\bf  R}^k$.
Put $X_i\equiv X_i({\bf z}):={I}(y_i\le\bf z)$, where ${I}(\cdot)$~
is the indicator function and the sign ``$\le$''
means the coordinatewise partial
order in $\bf  R^k$. We consider $X_i$
as elements of the Banach space
of all bounded right-continuous functions
(in the sense of the introduced partial order) defined on the set
$\{{\bf z}\in {\bf  R}^k: {\bf z}\le {\bf  a}\}$
and endowed with the sup-norm. Denote by ${\cal A}$
the cylindrical $\sigma$-field. Notice that the sup-norm is
${\cal A}$-measurable although the space ${\bf  B}$ is not separable. Without loss
of generality, we may assume that
${\bf  P}(X_i\ne0)\equiv{\bf  P}(y_i\le {\bf  a})\ne 0$.
Then the indicator random variables $X_i^0$   are well defined and meet
 (\ref{eqno12}). Therefore, from (\ref{eqno11}) and (\ref{eqno14}) we obtain the estimate
\begin{equation}
d (z,{S}_n,{\Pi}_n)\le\Delta(z,n,p),
                                                                \label{eqno15}
\end{equation}
where
$$
{S}_n\equiv{S}_n(t):= \lfloor nt\rfloor F_{\lfloor nt\rfloor }({\bf z}),
\quad
{\Pi}_n={\Pi}_n(t):= \sum\limits_{i=1}^{\overline\pi(\lfloor nt\rfloor )}{ I}(y_i\le{\bf z}),\,\,\,{\bf z}\le {\bf  a},\,t\in [0,1];
$$
$\Delta(z,n,p)$ is the right-hand side of (\ref{eqno14}), with
$p:={\bf  P}(y_1\le \bf  a)$.

The approximation of local empirical processes, in particular, in term of the distance $d (z,{S}_n,{\Pi}_n)$,
was studied by a number of authors (e.g., see Adell and de~la~Cal, 1994; Borisov, 1993, 1996, 2000; Horv{\'a}th, 1990; Major, 1990),
where relations similar to (\ref{eqno11}) were applied as well, which  allowed
to reduce the problem to the one-dimensional Bernoulli case.


The main result of the paper is the following generalization of Theorem~1, which allows us to use Lemma~3 in the non-i.i.d. case as well.

\pagebreak
{\bf Theorem 2.} {\it   For all integer $z\ge 0$,
\begin{equation} \label{nc:1}
d(z,\bar\nu_n(\cdot),\bar\pi_n(\cdot)) <
7\cdot 10^6 \left(p^*\right)^{z+1} \exp\left\{-\frac{1}{2}z\log\log(z+8)\right\}
\mbox{ if } \sum_{i=1}^n p_i^2 \le 1,
\end{equation}
where $p^*=\max\left\{\max_{i\le n}p_i,\sum_{i\le n}p_i^2\right\}$,
\begin{eqnarray} \label{nc:2}
d(z,\bar\nu_n(\cdot),\bar\pi_n(\cdot)) &<&
 10^8 \left( \sum_{i=1}^n p_i^2\right)^{(z+2)/2}
 \exp\left\{-\frac{1}{8}z\right\}\quad
 \mbox{ if } \sum_{i=1}^n p_i^2 \le 1,\\
 \label{nc:3}
d(z,\bar\nu_n(\cdot),\bar\pi_n(\cdot)) &<&
 14 \sum_{i=1}^n p_i^{z+2} \exp\left\{-\frac{1}{3} z \right\}\quad
\mbox{ if } \sum\limits_{i=1}^n p_i \le \frac{1}{2} .
\end{eqnarray}
Moreover,
\begin{equation} \label{nc:4}
d(0,\bar\nu_n(\cdot),\bar\pi_n(\cdot))\le \sum\limits_{i=1}^np_i^2.
\end{equation}
}

The following theorem is devoted to lower bounds for the distance $d(z,\bar\nu_n(\cdot),\bar\pi_n(\cdot))$.
This assertion generalizes Theorem 2 in Borisov (1996) where the i.i.d. case was considered.

{\bf Theorem 3.} {\it  For all integer $z\ge 0$,
\begin{equation} \label{nc:5}
d(z,\bar\nu_n(\cdot),\bar\pi_n(\cdot))\ge
\frac{1}{\sqrt {2\pi}}\exp\{-(z+3)\log(z+2)+z\} \sum_{k=1}^n B_k p_k^{z+2},
\end{equation}
where $B_k:=\exp\{-\sum_{i=1}^k p_i\}\prod_{j=1}^{k-1}p_j$, $k=1,\ldots,n$, and $\prod_{j=1}^{0}=1$ by definition.

Moreover,
\begin{equation} \label{nc:6}
d(0,\bar\nu_n(\cdot),\bar\pi_n(\cdot)) \ge 1-\exp \biggl\{-\frac{1}{2}\sum\limits_{i\le n}p_i^2(1-p_i)\biggr\}.
\end{equation}
}

Theorems 2 and 3 will be proved in Section~2.

R e m a r k. The lower bound (\ref{nc:5}) shows that, in general, the index $z+2$ on the right-hand side of (\ref{nc:3}) cannot be improved.
Moreover, observe that the  lower bound (\ref{nc:6}) has the asymptotic $\frac{1}{2}\sum\limits_{k\le n}p_k^2$ if $\sum\limits_{k\le n}p_k^2\to 0$.
In other words, the total variation distance $V(\bar\nu_n(\cdot),\bar\pi_n(\cdot))$ cannot decrease faster than $\sum\limits_{k\le n}p_k^2$ in contrast to the estimates in (\ref{eqno5}).

\newpage
\section{Proofs}

{\it Proof of Theorem }2.
We begin with constructing the families of random variables $\{\nu_i;\ i\le n\}$ and
$\{\pi_i;\ i\le n\}$ on a common probability space.
We will use  an approach, in fact, proposed
independently  by Serfling (1975) and Borovkov (1976, English translation: 2013). But to describe this approach we will prefer somewhat another terminology.
Given a random variable $\eta$ with distribution function $F_\eta(t)$, introduce the
quantile transform
$F_\eta^{-1}(\omega):=\inf\{t\in \overline {\bf   R}: F_\eta(t)\ge\omega\}$,
where $\overline {\bf   R}$~ is the extended real line
and $\omega\in[0,1]$. It is well known that if
$\omega$~ has the uniform distribution on $[0,1]$, then
${\cal L} (F_\eta^{-1}(\omega) )={\cal L}(\eta)$.
Now we set
$$
\nu_i^*:=F_{\nu_i}^{-1}(\omega_i),
\quad
\pi_i^*:=F_{\pi_i}^{-1}(\omega_i),
$$
where $\{\omega_i;\ i\le n\}$~are independent random variables with the uniform
distribution on $[0,1]$. One has
$$
\max\limits_{k\le n}
\left|\sum\limits_{i=1}^{k}\nu_i^*-\sum\limits_{i=1}^{k}\pi_i^*
\right|\le\sum\limits_{i=1}^{n}\zeta_i,
$$
where
$\zeta_i:=|\nu_i^*-\pi_i^*|=(\pi_i^*-1){\bf  I}\{\pi_i^*\ge 2\}
+{\bf  I} \{\omega_i\in[1-p_i,e^{-p_i}] \}$.
It is easy to verify that the distribution of $\zeta_i$ is determined as
follows:
$$
{\bf  P}(\zeta_i=0) =1-p_i(1-e^{-p_i}),\quad
{\bf  P}(\zeta_i=1) =e^{-p_i}-1+p_i+\frac{p_i^2}{2}e^{-p_i},
$$
\begin{equation*}
{\bf  P}(\zeta_i=k) =\frac{p_i^{k+1}}{(k+1)!}e^{-p_i},\quad k\ge 2.
\end{equation*}
Hence,
$${\bf  P}(\zeta_i=1)\le p_j^2,$$
\begin{equation}\label{zetai}
{\bf  P}(\zeta_i\ge k) < \frac{p_i^{k+1}}{(k+1)!\left(1-\frac{p_i}{k+2}\right)}e^{-p_i}
\le \frac{p_i^{k+1}}{(k+1)!}
,\quad k\ge 2.
\end{equation}
Besides, we will need the following inequality (see relation (28) in Borisov, 1996):
\begin{equation}\label{eet}
{\bf E} e^{t\zeta_i}\le \exp\left( e^t p_i^2 \exp(e^t p_i)\right)
\end{equation}
for $t>0$.

Prove now the relation (\ref{nc:1})
Denote $$t:=\log\left(\log(z+8)/2p^*\right).$$
Then, by \eqref{eet},
$$ d(z,\bar\nu_n(\cdot),\bar\pi_n(\cdot)) \le
\frac{{\bf E}\exp \left(t\sum_{i=1}^n \zeta_i\right)}{\exp(t(z+1))}\le
\exp\left\{e^t\sum_{i=1}^{n}p_i^2\exp\{e^tp_i\} - t(z+1) \right\} $$
$$ \le \exp\left\{\frac{1}{2}(z+8)^{1/2}\log(z+8)
-(z+1)\log\left(\log(z+8)/2p^*\right) \right\} $$
$$<7\cdot 10^6 \left(p^*\right)^{z+1} \exp\left\{-\frac{1}{2}z\log\log(z+8)\right\}, $$
that proves inequality (\ref{nc:1}).

The estimate (\ref{nc:4}) is immediate from the definition of random variables $\zeta_k$:
$$d(0,\bar\nu_n(\cdot),\bar\pi_n(\cdot))\le {\bf P}\Bigl (\bigcup_{k\le n} \{\zeta_k\neq 0)\}\Bigr)$$
$$ \le\sum_{k\le n}{\bf P} (\zeta_k\neq 0)=\sum\limits_{k\le n}p_k(1-e^{-p_k})
\le\sum\limits_{k\le n}p_k^2.
$$

Next, notice that  \eqref{nc:2} for $z=0$ follows from \eqref{nc:4}.
Thus, in order to prove \eqref{nc:2}, it remains to consider the case $z\ge 1$.

First, we will prove the estimate (\ref{nc:2}) by induction on $n$, under the condition
$$ \sum_{i\le n} p_i^2 \le \frac{1}{3} .$$
Denote $Q_n(k):={\bf P}\left(\sum_{i\le n} \zeta_i \ge k\right),\ k=1,2,\ldots$.
By the total probability formula one has
\begin{equation} \label{nn:0}
 Q_n(k) = {\bf P}(\zeta_n \ge k) + \sum_{m=0}^{k-1}
{\bf P}(\zeta_n=m)Q_{n-1}(k-m).
\end{equation}
It suffices to establish that, under the condition
 $ \sum_{i\le n} p_i^2 \le \frac{1}{3}$, the estimate
\begin{equation} \label{nn:1}
Q_n(k) < 5e^{-k/8}\left(\sum_{i=1}^n p_i^2\right)^{(k+1)/2}
\end{equation}
is valid for all integer $k \ge 2$.
For $n=1$ the relation (\ref{nn:1}) is immediate by \eqref{zetai}, since
$$
(k+1)!>4\, e^{k/8},\quad k\ge 2.
$$
Next, let (\ref{nn:1}) be valid for all $n \le N-1$.
We then obtain from (\ref{nn:0})
\begin{eqnarray*}
 Q_N(k) &
\le& Q_{N-1}(k)+\frac{p_N^{k+1}}{(k+1)!} +
p_N^2 Q_{N-1}(k-1)\\
&&+\frac{p_N^{k}}{k!} \sum_{i\le N-1} p_i^2+
\sum_{m=2}^{k-2} \frac{p_N^{m+1}}{(m+1)!} Q_{N-1}(k-m)\\
&=&
c_k \left( \sum_{i\le N}p_i^2\right)^{(k+1)/2}R,
\end{eqnarray*}
with
\begin{eqnarray*}
c_k&=& 5 e^{-k/8},\quad b=\frac{p_N^2}{\sum_{i\le N}p_i^2},\\
R&=&
 (1-b)^{(k+1)/2}+
\frac{1}{c_k(k+1)!} b^{(k+1)/2}+
\frac{c_{k-1}}{c_k} p_N b (1-b)^{(k-1)/2}\left( \sum_{i\le N-1} p_i^2 \right)^{1/2}
\\
&&+\frac{1}{c_k k!} b^{k/2}(1-b)^{1/2}\left( \sum_{i\le N-1} p_i^2 \right)^{1/2}
+\sum_{m=2}^{k-2}\frac{c_{k-m}}{c_k (m+1)!} b^{(m+1)/2} (1-b)^{(k-m)/2}
\\
&\le&
1-b+b \Bigg(
\frac{b^{(k-1)/2}}{c_k(k+1)!} + \frac{c_{k-1}}{c_k} \left( \sum_{i\le N-1} p_i^2 \right)^{1/2}
+\frac{b^{(k-2)/2}}{c_k k!} \left( \sum_{i\le N-1} p_i^2 \right)^{1/2}
\\
&&
\hspace{2cm} +\sum_{m=2}^{k-2}\frac{c_{k-m}b^{(m-1)/2}}{c_k (m+1)!}
\Bigg)
\\
&\le&
 1-b +b \Bigg(\frac{\sqrt{b}}{6c_2} + \frac{c_{k-1}}{c_k} 3^{-1/2}
+\frac{1}{4c_2} 3^{-1/2}+ \frac{c_{k-2}}{6c_k (1-c_{k-1}\sqrt{b}/(4 c_{k} ))} \Bigg)
\le 1,
\end{eqnarray*}
where we assumed, without loss of generality, that $b\le 1/2$ which can be achieved, e.g., by ordering
$p_1\ge p_2\ge p_3\ge \cdots$.

Thus, for $\sum_{i\le n}p_i^2\le \frac{1}{3}$ the estimate (\ref{nc:2}) is valid.
If $\frac{1}{3}\le \sum_{i\le n}p_i^2\le 1$ then, in order to prove the second estimate of the theorem, we use the first one:
$$ d(z,\bar\nu_n(\cdot),\bar\pi_n(\cdot))
<7\cdot 10^6 \exp\left\{-\frac{1}{2}z\log\log(z+8)\right\} $$
$$ < 10^8 \left(\sum_{i\le n}p_i^2\right)^{(z+2)/2}
\exp\left\{-\frac{1}{8}z\right\}. $$
The estimate (\ref{nc:2}) is proved.

Now, prove the estimate (\ref{nc:3}).
It suffices to establish that, under the condition $ \sum_{i\le n} p_i \le \frac{1}{2} $, the following inequality is valid:
\begin{equation} \label{nn:4}
Q_n(k) \le  19\sum_{i=1}^n (e^{-1/3}p_i)^{k+1}.
\end{equation}
Without loss of generality we can assume that $p_1\le p_2\le\ldots\le p_n$.
For $n=1$ the correctness of (\ref{nn:4}) is obvious.
Next, let (\ref{nn:4}) be valid for all $n \le N-1$. From  (\ref{nn:0}) we then obtain
\begin{eqnarray*}
Q_N(k) &\le& \frac{p_N^{k+1}}{(k+1)!} + 19\sum_{m=2}^{k-2}
\frac{p_N^{m+1}}{(m+1)!}\sum_{i\le N-1}(e^{-1/3}p_i)^{k-m+1} \\
&& + 19\sum_{i\le N-1}(e^{-1/3}p_i)^{k+1}
+19p_N^2 \sum_{i\le N-1}(e^{-1/3}p_i)^k
+ \frac{p_N^k}{k!}\sum_{i\le N-1}p_i^2 \\
& \le&
19 \sum_{i\le N}(e^{-1/3}p_i)^{k+1}R_1,
\end{eqnarray*}
where
\begin{eqnarray*}
 R_1 &=&
\frac{\sum_{i\le N-1}p_i^{k+1}}{\sum_{i\le N}p_i^{k+1}}+
\frac{e^{(k+1)/3}}{19(k+1)!}
\frac{p_N^{k+1}}{\sum_{i\le N}p_i^{k+1}} \\
&&+
\sum_{m=2}^{k-2} \frac{p_N^{m+1}
\sum_{i\le N-1}p_i^{k-m+1}}{\sum_{i\le N}p_i^{k+1}}\frac{e^{m/3}}{(m+1)!}
\\
&& + e^{1/3}p_N^2\frac{\sum_{i\le N-1}p_i^k}{\sum_{i\le N}p_i^{k+1}}
+ \frac{e^{(k+1)/3}p_N^k\sum_{i\le N-1}p_i^2}{19k!\sum_{i\le N}p_i^{k+1}}
\\
&\le& 1- \frac{p_N^{k+1}}{\sum_{i\le N}p_i^{k+1}}
\bigg( 1-\frac{e^{(k+1)/3}}{19(k+1)!}
- \frac{1}{2} \sum_{m=2}^{k-2}\frac{e^{m/3}}{(m+1)!}\\
&& \hspace{3cm} - \frac{1}{2}e^{1/3} - \frac{e^{(k+1)/3}}{2\cdot19k!} \bigg)
\le 1.
\end{eqnarray*}
Thus, the estimate (\ref{nc:3}) and Theorem 2 are proved.


\bigskip
{\it Proof of Theorem} 3.  Introduce the finite family of pairwise disjoint events
$$A_k:=\{\pi_1=1,\ldots,\pi_{k-1}=1, \pi_k=z+2\},\,\,\,k=1,\ldots,n.
$$
It is clear that for any construction of independent Bernoulli random variables $\{\nu_i\}$ on a common probability space with
independent Poissonian random variables $\{\pi_i\}$ the following implication of the events is valid:
$$\bigcup_{k\le n}A_k\subseteq \left\{\max_{k\le n}\left| \sum_{l\le k}(\nu_i-\pi_i)\right| >z\right\}.
$$
Hence,
$$d(z,\bar\nu_n(\cdot),\bar\pi_n(\cdot))\ge \sum_{k\le n}{\bf P}(A_k)=\frac{1}{(z+2)!}\sum_{i=1}^n B_k p_k^{z+2}.
$$
It remains to use Stirling's formula.

The last statement of the theorem is immediate from the simple lower bound
$$
d(0,\bar\nu_n(\cdot),\bar\pi_n(\cdot))\ge 1-{\bf  P} (\pi_i\in\{0,1\};\ i=1,\dots,n)=1-\prod\limits_{i=1}^ne^{-p_i}(1+p_i)$$
$$
\ge 1-\exp \biggl\{-\frac{1}{2}\sum\limits_{i\le n}p_i^2(1-p_i)\biggr\}$$
since
$$e^{-p_i}(1+p_i)\le \left(1-p_i+\frac{1}{2}p_i^2\right)(1+p_i)
\le 1-\frac{1}{2}p_i^2(1-p_i).$$

Theorem 3 is proved.

\begin{center}
{\bf References}
\end{center}



J.A.~Adell, J.~de~la~Cal (1996).
Optimal Poisson approximation of uniform empirical processes.
{\it Stochastic Processes and their Applications},
{\bf 64} (1), 135--142.

A.D.~Barbour, P.~Hall  (1984).
On the rate of Poisson convergence.
{\it Math. Proc. Cambridge Philos. Soc.}, {\bf 95}, 473--480.

A.D.~Barbour, L.~Holst, S.~Janson  (1992).
{\it Poisson Approximation}. Oxford Studies in Probability. 1st Edition. Clarendon Press.

I.S.~Borisov (1993).
Strong Poisson and mixed approximations of sums of
independent random variables in Banach spaces.
{\it Siberian Adv. Math.}, {\bf 3} (2), 1--13.

I.S. Borisov (1996). Poisson approximation of the partial sum
process in Banach spaces.
{\it Siberian Math. J.}, {\bf 37} (4), 627--634.

 I.S. Borisov (2000). A note on Poisson approximation of rescalled
set-indexed empirical processes. {\it Statist. Probab. Lett.} {\bf 46} (2), 101-103.

I.S.~Borisov, P.S.~Ruzankin (2002).
Poisson approximation for expectations of unbounded functions of independent random variables.
{\it Ann. Probab.}, {\bf 30} (4), 1657--1680.

I.S.~Borisov, I.S.~Vorozheikin (2008).
Accuracy of approximation in the Poisson theorem in terms of the $\chi^2$-distance.
{\it Siberian Math. J.}, {\bf 49} (1), 8--22.

A.A.~Borovkov (1976).
{\it Probability Theory}. Nauka, Moscow [in ~Russian].


A.A. Borovkov (2013). {\it Probability Theory}. Springer.


R.L.~Dobrushin (1970).
Prescribing a ~system of random variables
by conditional distributions,
{\it Theory Probab. Appl.},  {\bf 15} (3), 458--486.

R.M.~Dudley (1968). Distances of probability measures and random
variables. {\it Ann. Math. Statist.}, {\bf 39} (5), 1563--1572.

P.~Harremo\"{e}s, P.S.~Ruzankin (2004).
Rate of convergence to Poisson law in term of information divergence.
{\it IEEE Transac. Inform. Theory}, {\bf 50} (9), 2145--2149.

L.~Horv\'{a}th (1990).
A note on the rate of Poisson approximation of empirical processes,
{\it Ann. Probab.}, {\bf 18} (2), 724--726.

L.~Le Cam (1960).
An approximation theorem for the Poisson binomial distribution.
{\it Pacific~ J. Math.}, {\bf 10} (4), 1181--1197.

L.~Le Cam (1965).
On the distribution of sums of independent random variables.
in: {\it Bernoulli, Bayes, Laplace $($Anniversary Volume$)$},
Springer, 179--202.

L.~Le Cam (1970).
Remargues sur le th\'{e}or\`{e}me limit central dans les espaces
localement convexes.
in: {\it Les Probabilitit\'{e}s sur les Structures Alg\'{e}briques},  Paris,
C.N.R.S., 233--249.

P.~Major (1990).
A note of the approximation of the uniform empirical processes.
{\it Ann. Probab.}, {\bf 18} (1), 129--139.

S.Y.~Novak (2019a). On the accuracy of Poisson approximation. {\it Extremes}, epub ahead of print.\\
https://link.springer.com/article/10.1007/s10687-019-00350-6

S.Y.~Novak (2019b). {\it Poisson approximation}. Preprint. \\ https://arxiv.org/abs/1901.01847

Yu.V.~Prohorov (1953).
Asymptotic behavior of the binomial distribution.
{\it Uspekhi Mat. Nauk}, 8, 135--142 [in Russian].

P.S.~Ruzankin (2001).
On the Poisson Approximation of the Binomial Distribution.
{\it Siberian Math. J.}, {\bf 42} (2), 353--363.

P.S.~Ruzankin (2004).
On the Rate of Poisson Process Approximation to a Bernoulli Process.
{\it J. Appl. Probab.},
{\bf 41} (1), 271--276.

P.S.~Ruzankin (2010).
Approximation for expectations of unbounded functions of dependent integer-valued random variables.
{\it J. Appl. Probab.},
{\bf 47} (2), 594--600.

R.J.~Serfling  (1975).
A~general Poisson approximation theorem.
{\it Ann. Probab.}, {\bf 3} (3), 726--731.

V.M.~Zolotarev (1976).  Metric distances in spaces of random variables and their distributions. {\it Mathematics of the USSR -- Sbornik}, 30, No. 3, 373--401. \\doi:10.1070/sm1976v030n03abeh002280

\end{document}